\documentclass{ifacconf}  % Comment this line out if you need a4paper

\usepackage{mathtools}
\usepackage{amssymb}
\usepackage{floatrow}
\usepackage{xcolor}
\usepackage{graphicx}
\usepackage{natbib}    

\usepackage{mathrsfs} 
\usepackage{scalerel}
\newtheorem{lemma}{Lemma}
\newtheorem{dfn}{Definition}
\newtheorem{ass}{Assumption}

\DeclareMathOperator{\rank}{rank}

\newcommand\undermat[2]{%
  \makebox[0pt][l]{$\smash{\underbrace{\phantom{%
    \begin{matrix}#2\end{matrix}}}_{\text{$#1$}}}$}#2}

\newcommand{\RgeO}{\ensuremath{\mathbb{R}_{ \geq 0}}}
\newcommand{\Rat}[1]{\ensuremath{\mathbb{R}^{#1}}}

\begin{document}

\begin{frontmatter}

\title{State Estimation for a Class of Linear  Systems with Quadratic Output\thanksref{footnoteinfo}}

\thanks[footnoteinfo]{This research work is supported in part by NSERC-DG RGPIN-2020-04759, the SSF COIN project, the Swedish Research Council (VR), the Wallenberg AI, Autonomous Systems and Software Program
(WASP),  and the Knut och Alice Wallenberg foundation (KAW).}

\author[First]{Dionysios Theodosis} 
\author[Second]{Soulaimane Berkane} 
\author[Third]{Dimos V. Dimarogonas}

\address[First]{Dynamic Systems and Simulation Laboratory, School of Production Engineering and Management,
Technical University of Crete, Greece (e-mail: dtheodosis@dssl.tuc.gr).}
\address[Second]{D\'epartement   d'informatique   et   d'ing\'enierie, Universit\'e   du   Qu\'ebec   en   Outaouis,   Gatineau, Qu\'ebec,   Canada (e-mail: soulaimane.berkane@uqo.ca )}
\address[Third]{Division of Decision and Control Systems, School of Electrical Engineering and Computer Science, KTH Royal Institute of Technology, Sweden (e-mail: dimos@kth.se)}

\begin{abstract}
This paper deals with the problem of state estimation for  a class of linear time-invariant systems with quadratic output measurements. An immersion-type approach is presented that transforms the system into a state-affine system by adding a finite number of states to the original system. Under suitable persistence of excitation conditions on the input and its higher derivatives, global state estimation is exhibited by means of a Kalman-type observer. A numerical example is provided to illustrate the applicability of the proposed observer design for the problem of position and velocity estimation for a vehicle navigating in the $n-$dimensional Euclidean space using a single position range measurement.
\end{abstract}

\end{frontmatter}

\section{Introduction}
 
The aim of an observer is to provide an estimation of the running value of the system's internal
state using the input and output measurements. For linear systems, the observer synthesis is guaranteed through the observability property, namely, the determination of the initial state vector of the system from knowledge of the input and the corresponding output over an interval of time.  While observability is independent of the input for linear systems, this is not in general true for nonlinear systems and one needs to consider inputs that distinguish the states, namely,  inputs which generate different outputs, see \cite{HK77}.

Typically, the study of the observability of a nonlinear system is a local problem and can be characterized by the usual observability rank condition \cite{HK77}. However, this condition is not enough for the design of an observer since it tightly depends on the input. For such cases the design will be restricted to some appropriate classes of inputs, namely, regular or persistently exciting inputs, see for instance  \cite{B07}, \cite{BBH96}, \cite{BCC88}, \cite{GK01} and references therein. A well-known technique to design observers for nonlinear systems is the {\it immersion} approach where a nonlinear system is transformed into a state-affine system whose  dimension may be greater than the dimension of the initial system. Such methodologies have a long history. For instance,   \cite{FK83} presented a necessary and sufficient condition based on the observation space of the system. Another approach was considered in \cite{BS04} and \cite{J03} where the immersion was based on the solutions of a partial differential equation.  Another immersion-based technique was presented in \cite{BT07} for a wide class of (rank-observable) nonlinear
systems based on a high-gain design.

In this paper we consider systems with linear dynamics and   quadratic output measurements of the form \begin{align*}
&\dot{x}=Ax+Bu,\\
&y=\frac{1}{2}x^{\top}Cx,
\end{align*}
which is indeed a particular class of nonlinear systems. However, by restricting our attention to this class, our goal is to derive explicit conditions on the input $u$ that guarantee the design of an observer that is able to instantaneously estimate the state $x$ from the input and the (scalar) output measurement. First, through successive differentiation of the output, we extend the state of the system by a finite number of states which results in a new state-affine system with linear output. 
Then, we exploit the structure of the new extended system to derive suitable Persistence of Excitation (PE) conditions for the input and its derivatives that establish uniform observability for the new system. Consequently, the design of an observer for the obtained (uniformly observable) system follows directly from well known Kalman-like estimators, \cite{BBH96}, \cite{B07}, \cite{HS17} or other suitable observers. Since we consider an extended system, the estimate of the state of the original system can be obtained without any online inversion of a diffeomorphism. Finally, the framework presented in this paper generalizes and includes as a special case other state augmentation techniques presented in \cite{BSO11},   \cite{PI17}, \cite{HS17}, which mainly dealt with single and double integrator systems.  It should be noted that due to the nonlinear output of the considered class of systems it is also possible to apply other techniques as in  \cite{CDG93}, \cite{GK01},  \cite{GHO92} which for suitable inputs exploit a local change of coordinates to transform the system into a canonical form or by applying Lyapunov techniques as in \cite{T90}. {However, in contrast to these nonlinear techniques, our proposed approach has the advantage of employing a linear Kalman-type observer which guarantees global convergence while we also characterize explicitly the class of inputs (through the PE conditions) that guarantee the uniform observability property necessary for the exponential convergence of the estimator.} Finally, certain algebraic conditions for the observability of such systems was also proposed in \cite{D71}  without characterizing the admissible inputs. It should also be noted that the control of systems with quadratic outputs was considered in \cite{Monte17}.

\section{Preliminaries}\label{sec:prel}

\subsection{Notations} Throughout this paper we adopt the following notation. $\mathbb{N}$ and $\mathbb{R}$ denote, respectively, the sets of natural and real numbers.
For a given vector or matrix $(\cdot)\in{\mathbb R}^{n}$, $(\cdot)^{\top}$ denotes its transpose.
We denote by $I_{n}$ the $n\times n$ identity matrix. A matrix $A\in\Rat{n\times n}$ is called nilpotent if there exists an integer $\kappa\leq n$ such that $A^{\kappa}=0$.  By $0$ we denote each of the following: the scalar zero, the zero vector or the zero matrix. Depending on the context, the notation will be clear unless otherwise specified.  For $n\times n$ symmetric matrices $W$ and $Q$, the notation $W>0$ ($W\ge0$) is used if $W$ is positive definite (semi-definite) and $W>Q$ ($W\ge Q$) if $W-Q>0$ ($W-Q\ge 0$). With $\det(A)$ we denote the determinant of a square matrix  $A$. The Lie derivative of the real-valued scalar function $h$ along the vector field $f$ is denoted by $L_fh(x):=(dh(x)/dx)f(x)$ and the iterated derivatives  are defined as  $L^i_fh(x)=L_fL^{i-1}_fh(x)$, $L^0_fh(x)=h(x)$.

\subsection{Observability and Observers for LTV systems}

We first recall some well-known definitions and results for the observability of an LTV system.  Consider a linear time-varying system
\begin{subequations}\label{sys:LTV}
    \begin{align}
        \dot{x}=&A(t)x+B(t)u\label{sys:LTV:dynamics}\\
        y=& C(t)x\label{sys:LTV:out}
    \end{align}
\end{subequations}
where $x\in\Rat{n}$ is the state, $u\in\Rat{p}$ is the input, and $y\in\Rat{m}$ is the output of the system. $A(t), B(t), C(t)$ are matrix-valued functions of appropriate dimensions. We assume that these functions are continuous and bounded on $[0, +\infty)$.
\begin{dfn}\label{dfn:1}
The LTV system \eqref{sys:LTV} is called observable on $[t_0,T]$ if any initial state $x(t_0)=x_0$ is uniquely determined by the input $u(t)$ and the output $y(t)$ for $t\in[t_0,T]$.  $\triangleleft$
\end{dfn}

\begin{thm}\label{thm:1} [\cite{Rugh96}]
System \eqref{sys:LTV} is observable on $[t_0,T]$ if and only if the matrix:
$$W(t_0,T):=\int_{t_0}^{T}\Phi^{\top}(s,t)C^{\top}(s)C(s)\Phi(s,t)ds$$ is invertible, where
$\Phi(t,s)$ is the transition matrix defined by $\frac{d}{dt}\Phi(t,s)=A(t)\Phi(t,s)$, $\Phi(t,t)=I_n$. The matrix $W$ is called the Observability Gramian of \eqref{sys:LTV}.  $\triangleleft$
\end{thm}

Define the observability matrix
\begin{equation*}
\mathcal{R}(t):=\left(\begin{matrix}
N_0(t)\\N_1(t)\\\vdots\\N_{n-1}(t)
\end{matrix}\right)
\end{equation*}
with $N_0(t)=C(t)$, $N_{i}(t)=N_{i-1}(t)A(t)+\dot{N}_{i-1}(t)$, $i=1,\ldots,n-1$.
\begin{thm}\label{thm:2} [\cite{Rugh96}, \cite{SM67}]
The system \eqref{sys:LTV} is  observable if   $\rank \mathcal{R}(\bar{t})=n$, for some $\bar{t}\in [t_0,T]$.  $\triangleleft$
\end{thm}
The characterization of observability for time-varying systems is ``tied'' to finite time intervals, see \cite{BPP10}, \cite{Rugh96}, \cite{SM67}, and, \cite{W65} for different observability concepts and definitions.  For the state estimation problem a stronger property is required:
\begin{dfn}\label{dfn:2} [\cite{B07}] 
System \eqref{sys:LTV} or the pair $(A(t),C(t))$ is uniformly observable if there exist $\delta>0$, $\mu>0$ such that
\begin{equation}\label{obs:gramian}
\forall t\ge0\qquad W(t,t+\delta)\ge \mu I_n >0.  \triangleleft
\end{equation}
\end{dfn}

\begin{lemma}\label{lemma:Sca}
[\cite{Sc13}] Assume that there exists a positive integer $K$ such that the $k-$th derivative of $A$ (respectively $C$) is well defined and bounded up to $k=K$ (respectively up to $k=K+1$). If there exists a  matrix-valued function $\mathcal{O}(\cdot)$ of dimension $(\ell\times n)$, $\ell\ge1$, composed of row vectors of $N_1,\ldots,N_K$ such that for some strictly positive numbers $\bar{\delta},\bar{\mu}$ and $\forall t\ge0$ 
\begin{equation}\label{Sca}
 \int_t^{t+\bar{\delta}} \det( \mathcal{O}^{\top}(s)\mathcal{O}(s))ds\ge \bar{\mu}>0
\end{equation}
then system \eqref{sys:LTV} is uniformly observable. $\triangleleft$
\end{lemma}
 
A  Kalman-like observer for a uniformly observable LTV system \eqref{sys:LTV} is given in the following theorem:
\begin{thm}\label{Theorem:K:obs} [\cite{B07}]
If \eqref{sys:LTV} is uniformly  observable then there exists an observer of the form
$$\dot{\hat{x}}=A(t)\hat{x}+B(t)u+K(t)(y-C(t)\hat{x})$$
where
\begin{align*}
K(t)=&M(t)C^{\top}(t)W^{-1}\\
\dot{M}(t)=&A(t)M(t)+M(t)A^{\top}(t)\\&-M(t)C^{\top}(t)W^{-1}C(t)M(t) +V+\theta M(t)\\
M(0)=&M_0=M_0^{\top}>0,\,\,W=W^{\top}>0
\end{align*}
with $\theta>2||A(t)||$ for all $t\ge0$, or $V=V^{\top}>0$.  $\triangleleft$
\end{thm}

The boundedness assumption on $A(\cdot)$ and the uniform observability ensure that the solution $M(t)$ remains bounded for all times and the error $e:=\hat{x}-x$ between the state of the observer and the actual state decays exponentially to zero with the rate of convergence tuned by $\theta$ or $V$. For $\theta=0$ we obtain the usual Kalman observer; more details can be found in \cite{B07}.

\section{Problem formulation}
Consider the dynamical system
\begin{subequations}\label{sys:LTI}
    \begin{align}
        \dot{x}=&Ax+Bu\label{sys:LTI:dynamics}\\
        y=&\frac{1}{2}x^{\top}Cx \label{sys:LTI:out:q}
    \end{align}
\end{subequations}
where $x\in\mathbb{R}^n$ is the state, $u\in\Rat{p}$ is the input and $y\in\Rat{}$ is a scalar output. The constant matrices  $A\in\Rat{n\times n}$ and $B\in\Rat{n\times p}$ are arbitrary and, without loss of generality, the constant matrix $C\in\Rat{n\times n}$ is assumed symmetric. System \eqref{sys:LTI} is a linear time-invariant system with a single quadratic output and, thus, it is a special class of nonlinear systems. In contrast to classic linear systems, it is known that the observability of nonlinear systems depends usually on the input and is characterized locally. For instance, for the trivial system  $\dot{x}=0$, $y=x^2$, it is not possible to {\it distinguish} the initial conditions $x_0$ and $-x_0$ using only the output measurement. For certain nonzero input $u$, however, it is possible to distinguish all states of $\dot{x}=u$. Observability of nonlinear systems can be discussed using the notion of observation space where the observability rank condition can be used to study the so-called local weak observability around a given point, see \cite{HK77}, \cite{NS90}.

Notice that for zero inputs, the observation space is spanned by elements of the form $L_{f_0}^ih(x)$ with  $f_0(x)=Ax$ and $h(x)=\frac{1}{2}x^{\top}Cx$. In particular we have 
\begin{align}
	L_{f_0}h(x)=&x^{\top}CAx=\frac{1}{2}x^{\top}(CA+A^{\top}C)x:=\frac{1}{2}x^{\top}C_1x\nonumber\\
	L^2_{f_0}h(x)=&x^{\top}C_1Ax=\frac{1}{2}x^{\top}(C_1A+A^{\top}C_1)x:=\frac{1}{2}x^{\top}C_2x\nonumber\\
	&\vdots\nonumber\\
	L^i_{f_0}h(x)=&\frac{1}{2}x^{\top}C_{i}x,\quad i\in\mathbb{N}\label{Lie:der}
\end{align}
where the matrices $C_i$ are defined recursively as follows
\begin{align}
    &\left\{\begin{array}{l}
      C_0= C, \\
        C_{i+1}= C_iA+A^{\top}C_i.\\
    \end{array}\right. \label{Ci} 
\end{align}
Note that the matrix $C_i$, can be explicitly calculated using the following formula:
\begin{align}
    C_{i}=&\sum_{r=0}^i\binom{i}{r}{A^\top}^rCA^{i-r}\label{Ci:formula}
\end{align}
The proof of \eqref{Ci:formula} follows by simple induction; for completeness the proof can be found in the appendix.

In the subsequent sections, we will exploit the terms $L^{i}_{f_0}h(x)$   to augment the system \eqref{sys:LTI} with the additional states in order to bring the system in a new suitable form where an observer can be designed. Then, we will derive sufficient conditions for the admissible inputs that render the new extended time-varying system uniformly observable in the sense of Definition \ref{dfn:2}. To this end we start by our main assumption:

\begin{ass}\label{assumption:Cm}
  There exists $m\in\mathbb{N}$ with $C_{m}=0$.  
\end{ass}
Assumption \ref{assumption:Cm} is the only restriction we impose on the class of systems considered. The motivation behind this assumption is to facilitate the augmentation of the system by a finite number of states. Roughly speaking, this assumption is equivalent to the fact that there exists $m\in\mathbb{N}$ such that the $m-$th derivative $y^{(m)}$ of the output is zero under zero input or, equivalently,  $y^{(m)}$ is a polynomial function of time when $u=0$.  Note that for $u=0$, the solution of the linear time-invariant system is given by
\begin{align}
    x(t)=\exp(At)x(0),
\end{align}
which implies that 
\begin{align}
   y(t)=\frac{1}{2}x(0)^\top\exp(At)^{\top} C\exp(At)x(0).
\end{align}
For instance, it is clear that the output will be polynomial if the Taylor series defining the exponential matrix $\exp(At)$ is finite, i.e., when $A$ is nilpotent. This can also be seen  from \eqref{Ci:formula}, when $A^{k}=0$ then $C_{m}=0$ with $m=2k-1$. Also, when $A$ is skew-symmetric ($A=-A^{\top}$) and the matrices $A$, $C$ commute ($CA=AC$), one has $C_1=A^TC+CA=-CA+AC=0$, hence $m=1$ satisfies the assumption. The navigation example we provide in the simulation section also satisfies this assumption since the state matrix is nilpotent.
\section{State augmentation}
In this section, we proceed with the transformation of system  \eqref{sys:LTI} to an equivalent time-varying system when Assumption \ref{assumption:Cm} holds. We extend the state of the system   with $m$  additional states \begin{equation}\label{zi:state}
	z_i:=\frac{1}{2}x^{\top}C_ix  ,\,\,\,i=0,1,\cdots,(m-1).
\end{equation}
Then, since $C_{m}=0$, we have
\begin{align}\nonumber
\dot{z}_i=&\frac{1}{2}x^{\top}(A^{\top}C_i+C_iA)x+u^{\top}B^{\top}C_ix\\
=&z_{i+1}+ u^{\top}B^{\top}C_ix,\quad i=0,\ldots,(m-2),
\label{eq:dzi1}
\end{align}
and
\begin{equation}\label{eq:dzi2}
\dot{z}_{m-1}=u^{\top}B^{\top}C_{m-1}x,
\end{equation}
where the last equation holds due to Assumption \ref{assumption:Cm}.
Therefore, by defining the extended state as
\begin{equation}\label{ext:state:1}
z:=\left[\begin{matrix}
z_0&z_1& \cdots&z_{m-1}&x
\end{matrix}\right]^{\top}\in\mathbb{R}^{m+n}
\end{equation}
and in view of \eqref{sys:LTI} and \eqref{eq:dzi1}-\eqref{eq:dzi2}, the dynamics of the new variable $z$ are given by the following LTV system
\begin{subequations}\label{sys:noC:ext}
\begin{align}
\dot{z}=&\mathcal{A}(u)z+\mathcal{B}u\label{sys:noC:ext:s},\\
y=&\mathcal{C}z \label{sys:noC:ext:out},
\end{align}
\end{subequations}
where the matrices $\mathcal{A}(u)\in\mathbb{R}^{(m+n)\times(m+n)}$, $\mathcal{B}\in\mathbb{R}^{(m+n)\times p}$ and $\mathcal{C}\in\mathbb{R}^{1\times(m+n)}$ are given by
\begin{align}
\mathcal{A}(u):=&\left[\begin{matrix}
0&1&0&\cdots&0&u^{\top}B^{\top}C_0\\
0&0&1&\cdots&0&u^{\top}B^{\top}C_1\\
\vdots&\vdots&\vdots&\vdots&\vdots\\
0&0&0&\cdots&1&u^{\top}B^{\top}C_{m-2}\\
0&0&0&\cdots&0 &u^{\top}B^{\top}C_{m-1}\\
0_{n\times 1}&0_{n\times 1}&0_{n\times 1}&0_{n\times 1}&0_{n\times 1}&A
\end{matrix}\right]\label{matrix:Au:ext}\\
\mathcal{B}:=&\left[\begin{matrix}
0_{m\times p}\\B
\end{matrix}\right] \label{matrix:B:ext}\\ 
\mathcal{C}:=&\left[\begin{matrix}  \begin{matrix}
1&0&0&\cdots&0
\end{matrix}
&0_{1\times n}\end{matrix}\right]\label{matrix:Cu:ext}.
\end{align}
Notice that the new augmented system \eqref{sys:noC:ext:s} is a state affine system, \cite{BBH96}, which can also be considered as a LTV system for some fixed input function $u(t)$. We adopt the state-affine definition to emphasize the dependence on the input $u(\cdot)$ even if $u$ appears linearly in $\mathcal{A}$. 
For state affine systems several Kalman-type observer designs have appeared in the literature, see for instance   \cite{BCC88}, \cite{B07}, \cite{BBH96}. Typically, the main property required to  use a Kalman-type observer is that of uniform observability. Therefore, to estimate the state of the extended system, it suffices to study the observability of the pair $(\mathcal{A}(u),\mathcal{C})$, see Section \ref{section:observability}.

\section{Uniform observability analysis}\label{section:observability}
In this section, we derive different sufficient conditions for the admissible inputs $u(\cdot)$ that render the system \eqref{sys:noC:ext} uniformly observable. Before we proceed, the following necessary condition is provided, which allows the input $u$ to directly affect the extended system and its observability properties.  
\begin{prop}\label{prop:consistency}
If $u(t)\equiv 0$ for all times or $B^{\top}C_r=0$ for all $r\in\mathbb{N}$, then the pair $(\mathcal{A}(u),\mathcal{C})$ is not uniformly observable. 
\end{prop}

\begin{pf}
Notice that for the cases $u\equiv0$ or $B^{\top}C_r=0$ for all $r\in\mathbb{N}$, the matrix $\mathcal{A}(u(t))=\mathcal{A}$ is constant. In that case, the pair $(\mathcal{A},\mathcal{C})$ in \eqref{matrix:Au:ext}, \eqref{matrix:Cu:ext} is not Kalman observable since the observability rank condition gives $\rank(\mathcal{C}^\top, \mathcal{A}^\top\mathcal{C}^\top,\cdots,\mathcal{A}^{{m+n-1}^\top}\mathcal{C}^\top)=m$. To prove this claim rewrite $\mathcal{A}$ in the following block structure
\begin{equation*} 
\mathcal{A} =\begin{bmatrix}
S_m&0\\0&A
\end{bmatrix}
\end{equation*}
where $S_m\in\Rat{m\times m}$ is the standard shift matrix 
\begin{align}\label{mat:shift}
S_m= \begin{bmatrix}
0&1&0&\ldots &0\\
0&0&1&\ldots &0\\
\vdots&\vdots&\vdots&\ddots &0\\
0&0&0&\ldots &1\\
0&\undermat{m-1}{0&0&\ldots &0}
\end{bmatrix}\\ \nonumber
\end{align}
Notice now that due to the triangular block structure we also have that
$$\mathcal{A}^j=\begin{bmatrix}
S_m^j&0\\0&A^j
\end{bmatrix}$$
Then, it follows by direct calculations that $\mathcal{C}\mathcal{A}=(0,1,0,\ldots,0)$, $\mathcal{C}\mathcal{A}^2=(0,0,1,\ldots,0)$, $\ldots,   \mathcal{C}\mathcal{A}^{m-1}=(0,0, \ldots,1,0)$ and hence $\rank(\mathcal{C}^\top, \mathcal{A}^\top\mathcal{C}^\top,\cdots,\mathcal{A}^{{m-1}^\top}\mathcal{C}^\top)=m$. Since $S_m$ is a shift matrix we also have that $S_m^m=0$. The latter in conjunction with \eqref{matrix:Cu:ext} implies that $\mathcal{C}\mathcal{A}^m=0$ and $\mathcal{C}\mathcal{A}^j=0$ for $j=m+1,\ldots,m+n-1$ which proves that the system is not rank observable.
\end{pf}

The state-transition matrix associated with \eqref{sys:noC:ext} is defined by \begin{equation}\label{mat:trans}
    \frac{d}{dt}\Phi(t,\tau)=\mathcal{A}(u(t))\Phi(t,\tau), \quad \Phi(\tau,\tau)=I_{m+n}.
\end{equation}  
In general, calculating the transition matrix $\Phi(t,\tau)$ and verifying that the inequality $W(t,t+\delta)\ge \mu I$  holds is a tedious task, especially in our case where the state matrix $\mathcal{A}(u(t))$ depends on the input. However, we can exploit the block structure of \eqref{matrix:Au:ext} and simplify its representation. More specifically, rewrite $\mathcal{A}(u)$ in the following form
\begin{equation}\label{mat:Au:small}
\mathcal{A}(u(t))=\begin{bmatrix}
S_m&U(t)\\0&A
\end{bmatrix}
\end{equation}
where $S_m\in\Rat{m\times m}$ is  given by \eqref{mat:shift}
and $U:\RgeO\to\Rat{m\times n}$ is given by
\begin{equation}\label{mat:U}
U(t)=\begin{bmatrix}
u^{\top}(t)B^{\top}C_0\\
u^{\top}(t)B^{\top}C_1\\
\vdots\\
u^{\top}(t)B^{\top}C_{m-2}\\
u^{\top}(t)B^{\top}C_{m-1}\\
\end{bmatrix}.
\end{equation}
Notice that due to the structure of $\mathcal{A}(u(t))$ in \eqref{mat:Au:small} it follows from \eqref{mat:trans}, that the transition matrix has the following form:
\begin{equation}\label{transition:mat}
    \Phi(t,\tau)=\begin{bmatrix}
\Phi_{11}(t,\tau)&\Phi_{12}(t,\tau)\\
0&\Phi_{22}(t,\tau)
\end{bmatrix}
\end{equation}
with
\begin{align}
&\frac{d}{dt}\Phi_{11}(t,\tau)=S_m\Phi_{11}(t,\tau), \quad\Phi_{11}(\tau,\tau)=I_{m},\label{Phi:11:de}\\
&\frac{d}{dt}\Phi_{22}(t,\tau)=A\Phi_{22}(t,\tau), \quad \Phi_{22}(\tau,\tau)=I_{n},\label{Phi:22:de}
	\end{align}  
and $\Phi_{12}(t,\tau)$ satisfying
\begin{equation}\label{Phi:12:de}
    \frac{d}{dt}\Phi_{12}(t,\tau)=S_m\Phi_{12}(t,\tau)+U(t)\Phi_{22}(t,\tau),\quad \Phi_{12}(\tau,\tau)=0.
\end{equation}
Since $S_m$ and $A$ are constant matrices, we have that
\begin{subequations}
	\begin{align}
	\Phi_{11}(t,\tau)=&\exp\{S_m(t-\tau)\},\label{Phi:11}\\
	\quad\Phi_{22}(t,\tau)=&\exp\{A(t-\tau)\}, \label{Phi:22}
	\end{align}
whereas $\Phi_{12}(\cdot,\cdot)$ has the following form
\begin{equation}\label{Phi:12}
\Phi_{12}(t,\tau)=\int_{\tau}^t\Phi_{11}(t,s)U(s)\Phi_{22}(s,\tau)ds.    
\end{equation}
\end{subequations}
Indeed, $\Phi_{12}(\tau,\tau)=0$ and from \eqref{Phi:11:de} and the Leibniz integral rule we also have 
\begin{align*}
    \frac{d}{dt}\Phi_{12}(t,\tau)=&\Phi_{11}(t,t)U(t)\Phi_{22}(t,\tau)\\&+\int_{\tau}^{t}\frac{d}{dt}\Phi_{11}(t,s)U(s)\Phi_{22}(s,\tau)ds\\
    =&S_m\int_{\tau}^{t}\Phi_{11}(t,s)U(s)\Phi_{22}(s,\tau)ds\\&+U(t)\Phi_{22}(t,\tau)\\
    =&S_m\Phi_{12}(t,\tau)+U(t)\Phi_{22}(t,\tau),
\end{align*}
so that \eqref{Phi:12:de} follows as well.
By defining $\mathcal{C}=(\mathcal{C}_m,0_{1\times n})$, with $\mathcal{C}_m=(1,0_{1\times m-1})$ and taking into account the observability condition in Definition \ref{dfn:2} and \eqref{matrix:Cu:ext}, we have that the Observability Gramian of the extended system \eqref{sys:noC:ext} is given by
	\begin{align}
	&W(t,t+\delta)=\nonumber\\
	&\int_t^{t+\delta}\begin{bmatrix}
	\Phi_{11}(s,t)^{\top}\\
	\Phi_{12}(s,t)^{\top}
	\end{bmatrix}
	\mathcal{C}_m^{\top}\mathcal{C}_m
	\begin{bmatrix}
	\Phi_{11}(s,t)&\Phi_{12}(s,t)
	\end{bmatrix}ds.
\label{obs:Gramian:ext}
	\end{align}
Notice that the Observability Gramian is expressed in terms of the matrices $S_m$, $A$ and $U(t)$ through the definitions of $\Phi_{11}(t,\tau)$, $\Phi_{22}(t,\tau)$ and $\Phi_{12}(t,\tau)$ above and does not require the derivatives of $u$ neither the evaluation of the usual Peano-Baker series, see \cite{Rugh96}. According to Theorem \ref{Theorem:K:obs}, to  design an observer for the time-varying system \eqref{sys:noC:ext}, we require   uniform observability, which guarantees the exponential convergence of the observer. Therefore, according to Definition \ref{dfn:2} it suffices to consider inputs $u$ for which the Observability Gramian \eqref{obs:Gramian:ext} of the extended system satisfies inequality \eqref{obs:gramian}, i.e., persistently exciting inputs, see \cite{B07}. The transition matrices $\Phi_{11}$ and $\Phi_{22}$ can be easily  computed from \eqref{Phi:11} and \eqref{Phi:22} since $S_m$ and $A$ are constant, however, verifying that the Observability Gramian in \eqref{obs:Gramian:ext}  satisfies inequality \eqref{obs:gramian} is  non-trivial.

Typically, to design a Kalman-type observer and guarantee its exponential convergence it is required that the control input $u(t)$ is bounded for all $t\ge0$ which also implies that $\mathcal{A}(u(t))$ is bounded, see Theorem \ref{Theorem:K:obs}. To derive more explicit uniform observability conditions, we further assume in this work that the higher derivatives of the inputs are bounded.
\begin{ass}\label{As:bounded:input}
	The input $u(t)$ is bounded and there exists  $\kappa\ge m$ such that the  derivatives $\dot{u}(t)$, $\ddot{u}(t),\cdots,  \overset{(\kappa)}{u}(t)$, are also bounded for all $t\ge0$.
\end{ass}

{The next result is based on Lemma \ref{lemma:Sca} and gives a sufficient Persistence of Excitation (PE) condition for the uniform observability of the extended system by exploiting row vectors of the observability matrix at the expense of sufficiently smooth input $u(\cdot)$ with bounded derivatives.}
\begin{prop}\label{PE:cond:2} 
	Let $r_0=0$ and
		\begin{equation}\label{fun:r:2}
		r_{i+1}(t)=r_i(t)A+\dot{r}_i(t)+u^{\top}(t)B^{\top}C_{i},\,\,i=0,1,\ldots.
		\end{equation} 
	and assume also that there exist positive constants $\delta$, $\mu$, $\kappa$ such that for all $t\ge0$ Assumption \ref{As:bounded:input} holds and the following condition is satisfied:
	\begin{align}
	\int_t^{t+\delta}\det\left(\sum_{i=m}^{\kappa}r_{i}^{\top}(s) r_{i}(s)\right)ds\ge\mu.	\label{uni:obs:2}
	\end{align}
	Then, the system \eqref{sys:noC:ext} is uniformly observable. 
\end{prop} 
	
\begin{pf}
	To show that \eqref{uni:obs:2} implies uniform observability of the system \eqref{sys:noC:ext}, we will exploit Lemma \ref{lemma:Sca}. Hence, as in the statement of Lemma \ref{lemma:Sca} we define the row vectors
	\begin{align}
    &N_0(t):=\mathcal{C}\nonumber\\
	&N_{i}(t):=N_{i-1}(t)\mathcal{A}(u(t))+\dot{N}_{i-1}(t),\,\,i=1,\ldots,\label{Ni}
	\end{align}
	In view of \eqref{matrix:Au:ext}-\eqref{matrix:Cu:ext} we obtain
	\begin{align*}
	N_0(t)=&\left[\begin{matrix}1&0&0&\cdots&0&0_{1\times n}\end{matrix}\right]\\
	N_1(t)=&\left[\begin{matrix}0&1&0&\cdots&0&r_1(t)\end{matrix}\right]\\
	\vdots&\\
	N_{m-1}(t)=&\left[\begin{matrix}0&0&0&\cdots&1&  r_{m-1}(t)\end{matrix}\right] \\
	N_{\kappa}(t)=&\left[\begin{matrix}0&0&0&\cdots&0&  r_{\kappa}(t)\end{matrix}\right],\,\,\,\kappa\ge m 
	\end{align*}
	where  $r_i(\cdot)$, $i=1,\ldots,\kappa$ is a sequence of $1\times n$ vectors defined by \eqref{fun:r:2}.
	Now, consider the matrix $\mathcal{O}(t)=\left[\begin{matrix}
	N_0(t)^\top&N_1(t)^\top&\ldots&N_{\kappa}(t)^\top\end{matrix}\right]^{\top}$
	from which we obtain the $(m+n)\times(m+n)$ matrix
	\begin{align*}
	& \mathcal{O}(t)^{\top}\mathcal{O}(t)= \\
	&\left[
	\begin{matrix}
	1&0&0&\cdots&0 &0_{1\times n}\\
	0&1&0&\cdots&0 &r_1(t)\\
	0&0&1&\cdots&0 &r_2(t)\\
	\vdots&\vdots&\vdots&\cdots& \vdots&\vdots\\
	0&0&0&\cdots&1&r_{m-1}(t)\\
	0_{n\times 1}&r_1^{\top}(t)&r_2^{\top}(t)&\cdots &r_{m-1}^{\top}(t) &\sum_{i=1}^{\kappa}r_i^{\top}(t)r_i(t)
	\end{matrix}
	\right]
	\end{align*}
	By taking into account Assumption \ref{As:bounded:input}, it follows that $\mathcal{A}(u(\cdot))$ and its $\kappa$- derivatives are bounded and thus all assumptions of Lemma \ref{lemma:Sca} hold. Hence, to show the uniform observability of the system it suffices to show that \eqref{Sca} holds. Notice that for the matrix $\mathcal{O}(\cdot)^{\top}\mathcal{O}(\cdot)$ above we have by the Schur complement and its determinant  that 
	\begin{align*}
	\det(\mathcal{O}^{\top}(t)\mathcal{O}(t))=& \det\left(\sum_{i=1}^{\kappa}r_i^{\top}(t)r_i(t)-\sum_{i=1}^{m-1}r_i^{\top}(t)r_i(t)\right)\\
	&=\det\left(\sum_{i=m}^{\kappa}r_i^{\top}(t)r_i(t)\right)
	\end{align*}
	Therefore, it follows that if condition \eqref{uni:obs:2} holds then also \eqref{Sca} is satisfied and from Lemma \ref{lemma:Sca} we conclude that the system \eqref{sys:noC:ext} is  uniformly observable. \qed
\end{pf}

{Notice that condition \eqref{uni:obs:2} requires a sufficiently smooth and bounded input as there is no restriction on how large the constant $\kappa$ may be. In particular, $\kappa$ must be at least greater or equal to $m+1$ since $\det(r_m^\top(s) r_m(s))=0$.  In the following material, we will derive PE conditions that require less number of derivatives. The first proposition below provides a PE condition for the admissible inputs and their derivatives and is equivalent to the uniform observability of the pair $(A,r_m(t))$, where $r_m(\cdot)$ is defined in \eqref{fun:r:2}.}
\begin{prop}\label{PE:cond:3} 
	Suppose that Assumption \ref{assumption:Cm} holds and in addition Assumption \ref{As:bounded:input} is satisfied with $\kappa=m$. Then, the system \eqref{sys:noC:ext} is uniformly observable if there exist positive constants $\bar{\delta}$ and $\bar{\mu}$ such that for all $t\ge0$ we have
		\begin{align}
		\int_t^{t+\bar{\delta}}\Phi_{22}^{\top}(s,t)r_m^{\top}(s) r_m(s)\Phi_{22}(s,t)ds\ge\bar{\mu} I_n.	\label{uni:obs:3}
		\end{align}
\end{prop} 

\begin{pf}
To show that the system \eqref{sys:noC:ext} is uniformly observable, it suffices to show that there exist $\delta$, $\mu>0$ such that \eqref{obs:gramian} holds, i.e., $W(t,t+\delta)>\mu I_{m+n}$, for all $t\ge0$ with $W(t,t+\delta)$ given by \eqref{obs:Gramian:ext}.
We proceed by contradiction. Suppose that for every $\mu>0$ and $\delta>0$ there exists $t\ge0$ with $W(t,t+\delta)<\mu I_{m+n}$. Consider a sequence $\{\mu_p\}_{p\in\mathbb{N}}$ that converges to zero  with $\mu_p>0$ and $\delta=\bar{\delta}>0$, with $\bar{\delta}$ satisfying the PE condition \eqref{uni:obs:3}.  Then, there exists a sequence of times $\{t_p\}_{p\in\mathbb{N}}$ and a sequence $\{\hat{d}_p\}_{p\in\mathbb{N}}$ with $\hat{d}_p\in D=\{d\in\mathbb{R}^{m+n}:||{d}||=1\}$ such that, for all $p\in\mathbb{N}$, $\hat{d}_p^\top W(t_p,t_p+\bar{\delta})\hat{d}_p<\mu_p$. Next, consider a sub-sequence of $\{\hat{d}_p\}_{p\in\mathbb{N}}$ which converges to some $d\in D$, since $D$ is compact. Let
$d=(d_1^{\top},d_2^{\top})^{\top}\in\Rat{m+n}$ with $d_1=(d_{11},d_{12},\ldots,d_{1m})^{\top}\in\Rat{m}$, $d_2=(d_{21},d_{22},\ldots,d_{1n})^{\top}\in\Rat{n}$ and $||d||=1$, then it follows from the previous assumption, \eqref{transition:mat}, and \eqref{obs:Gramian:ext} that
\begin{align}
    \lim_{p\to+\infty}\int_{t_p}^{t_p+\bar{\delta}}\|\mathcal{C}\Phi(s,t_p)d\|^2ds=0
\end{align}
or, by a change of variables,
\begin{align}
    \lim_{p\to+\infty}\int_{0}^{\bar{\delta}}\|f_p(s)\|^2ds=0
    \end{align}
where we define
\begin{align*}
    f_p(t)&=\mathcal{C}\Phi(t+t_p,t_p)d\\
          &=\mathcal{C}_m\Phi_{11}(t+t_p,t_p)d_1+\mathcal{C}_m\Phi_{12}(t+t_p,t_p)d_2,
\end{align*}
with $\Phi_{11}(\cdot,\cdot)$ and  $\Phi_{12}(\cdot,\cdot)$ satisfying \eqref{Phi:11} and \eqref{Phi:12}, respectively.
From \eqref{Phi:11:de}-\eqref{Phi:12:de}, the successive time-derivatives of $f_p(t)$ are given as follows (argument for $(t+t_p,t_p)$ in $\Phi_{ij}$ is omitted):
\begin{align*}
    f_p^{(1)}(t)&=\mathcal{C}_mS_m\Phi_{11}d_1+\mathcal{C}_mS_m\Phi_{12}d_2+r_1\Phi_{22}d_2\\
    f_p^{(2)}(t)&=\mathcal{C}_mS_m^2\Phi_{11}d_1+\mathcal{C}_mS_m^2\Phi_{12}d_2+r_2\Phi_{22}d_2\\
    &\vdots \\
    f_p^{(m)}(t)&=\mathcal{C}_mS_m^m\Phi_{11}d_1+\mathcal{C}_mS_m^m\Phi_{12}d_2+r_m\Phi_{22}d_2\quad\\
    &=r_m\Phi_{22}d_2\qquad (S_m^m=0)
\end{align*}
with $r_i$ defined by \eqref{fun:r:2}.
Now, using the result of Lemma A.1 in \cite{Sc13}, we deduce that 
\begin{align}\label{eq:int-fk}
    \lim_{p\to+\infty}\int_{0}^{\bar{\delta}}\|f_p^{(k)}(s)\|^2ds=0,\quad k=0,\cdots, m
    \end{align}
or, in particular for $k=m$,
\begin{align}
    \lim_{p\to+\infty}\int_{t_p}^{t_p+\bar{\delta}}\|r_m(s)\Phi_{22}(s,t_p)d_2\|^2ds=0.
\end{align}
However thanks to the PE condition \eqref{uni:obs:3}, this cannot hold except with $d_2=0$. This implies that the derivatives above have the form:
\begin{align}
    f_p^{(k)}(t)=\mathcal{C}_mS_m^k\Phi_{11}(t+t_p,t_p)d_1,\quad k=0,\cdots,m.\label{f:k:der}
\end{align}
On the other hand, notice that thanks to the special structure of the matrix $S_m$, the state transition matrix $\Phi_{11}(t,\tau)$ can be written as 
\begin{align*}
    \Phi_{11}(t,\tau)&=\exp(S_m(t-\tau))\\
    &=\begin{bmatrix}
    1&(t-\tau)&\frac{(t-\tau)^2}{2!}&\cdots&\frac{(t-\tau)^{m-1}}{(m-1)!}\\
    0&1&(t-\tau)&\cdots&\frac{(t-\tau)^{m-2}}{(m-2)!}\\
    \vdots&\vdots&\ddots&\vdots&\vdots\\
    0&0&0&\cdots&1\\
    \end{bmatrix}.
\end{align*}
Also, since $S_m$ is a shift matrix, we have 
\begin{align*}
    \mathcal{C}_mS_m&=[0\; 1\; 0\cdots0\; 0]\\
    \mathcal{C}_mS_m^2&=[0\; 0 \;1\cdots0\; 0]\\
      &\vdots\\
       \mathcal{C}_mS_m^{m-1}&=[0\; 0\;\ldots\;0 \;1].
\end{align*}
By taking into account \eqref{f:k:der}, the last equation $\mathcal{C}_mS_m^{m-1}=[0,0,\ldots,0,1]$ implies that $f_p^{(m-1)}(s)=d_{1,m}$ which in turn leads in view of \eqref{eq:int-fk} to $\lim_{p\to\infty}\int_{0}^{\bar{\delta}} |f_p^{(m-1)}(s)|^2ds=\lim_{p\to\infty}\allowbreak\int_0^{\bar{\delta}} d_{1,m}^2ds=0$ and $d_{1,m}=0$. Since $d_{1,m}=0$, we have in view of \eqref{f:k:der} and $\mathcal{C}_mS_m^{m-2}=[0,0,\ldots,1,0]$  that $f^{(m-2)}=(t-\tau)d_{1,m}+d_{1,m-1}=d_{1,m-1}$. This leads in view of \eqref{eq:int-fk} to  $\lim_{p\to\infty}\int_0^{\bar{\delta}}\allowbreak |f_p^{(m-2)}(s)|^2ds=\lim_{p\to\infty}\allowbreak \int_{0}^{\bar{\delta}}d_{1,m-1}^2ds=0$ and hence $d_{1,m-1}=0$. Using recursively the same argument for each $k=m-2,m-3,\ldots,1,0$ and by exploiting \eqref{f:k:der} we have that $d_1=0$ and thus $d=0$ which is a contradiction since $||d||=1$. \qed
\end{pf}

{Notice that the PE condition in \eqref{uni:obs:3} is equivalent to the uniform observability of the pair $(A,r_m(t))$ since the matrix in \eqref{uni:obs:3} corresponds to the observability Gramian of this pair (recall that $\Phi_{22}$ is the state transition matrix of $A$). The next PE condition guarantees the uniform observability of system \eqref{sys:noC:ext} without requiring the computation of the Gramian matrix for the pair $(A,r_m(t))$ but under the additional assumption that the matrix $A$ has real eigenvalues.}
\begin{prop}\label{PE:cond:4} 
When $A$ has real eigenvalues, a sufficient condition for \eqref{uni:obs:3} to hold is the following
\begin{align}\label{uni:obs:4}
		\int_t^{t+\bar{\delta}}r_m^{\top}(s) r_m(s)ds\ge\bar{\mu} I_n.
\end{align}
for some $\bar{\delta},\bar{\mu}>0$.
\end{prop} 
\begin{pf}
The proof follows directly from Lemma 2.7 in \cite{HS17}  by noticing that the pair $(A,I_n)$ is Kalman observable. \qed
\end{pf}

Notice that each $r_i(\cdot)$ in \eqref{fun:r:2}, can be written as a summation of higher order derivatives of the input $u$. More specifically, it can shown by induction that
\begin{align}
    r_m=\sum_{i=0}^{m-1}u^{(i)^\top}\Gamma_{i+1,m}
\end{align}
such that 
\begin{align}
    \Gamma_{i+1,k+1}=\Gamma_{i+1,k}A+\Gamma_{i,k},\quad i=0,\cdots,k
\end{align}
and $\Gamma_{0,k}=B^\top C_k$ and $\Gamma_{k+1,k}=0$.  It follows that 
\begin{align}\label{eq:rm:Gamma}
    r_m^\top=\begin{bmatrix}
    \Gamma_{1,m}\\
    \Gamma_{2,m}\\
    \vdots\\
    \Gamma_{m,m}
    \end{bmatrix}^\top\begin{bmatrix}
    u\\
    \dot u\\
    \vdots\\
    u^{(m-1)}
    \end{bmatrix}=:\Gamma^\top\bar U
\end{align}

\begin{prop}
Assume that $A$ has real eigenvalues and the pair $(A,\Gamma)$ is Kalman observable. Then, a sufficient condition for \eqref{uni:obs:3} to hold is the following
\begin{align}
		\int_t^{t+\delta}\bar U(s)\bar U(s)^\top ds\ge\mu I_{pm}.\label{PE:Ubar}
\end{align}
then \eqref{uni:obs:4} holds. 
\end{prop}
\begin{pf}
The proof follows directly from Lemma 2.7 in \cite{HS17} in view of \eqref{uni:obs:3} and \eqref{eq:rm:Gamma}. \qed
\end{pf}
The advantage of the above PE condition is that the vector $\bar U(t)$ is directly expressed as a function of the input and its higher derivatives. However, the drawback compared to \eqref{uni:obs:4} is that we might check a PE condition on a vector with higher dimension when $pm\geq n$. Finally, notice that conditions \eqref{uni:obs:3},  \eqref{uni:obs:4}, and \eqref{PE:Ubar} require in general less number of derivatives compared to \eqref{uni:obs:2}.

{Now that we established different conditions for the uniform observability of the extended system \eqref{sys:noC:ext} it is possible to estimate its state and consequently the state of the original system with the desing of a Kalman estimator as in Theorem \ref{Theorem:K:obs} with $\theta=0$.  Note that the structure of the augmented system allows to apply other observer designs presented for instance in \cite{BBH96}, \cite{BCC88}, \cite{KJ11}, \cite{TK19}. In fact, one of the advantages of the observer design methodology adopted in this paper is that the global state estimation can be achieved by a simple linear Kalman type observer as in Theorem \ref{Theorem:K:obs} as shown in the simulation examples of the next section.}

 \section{Numerical Example}
 The following example illustrates the state extension as well as the observability conditions presented in the previous section. We consider a vehicle navigating in $\mathbb{R}^n$ using a single position range measurement positioned at $0\in\Rat{n}$. The dynamics of the vehicle can be written as 
 \begin{align}
     &\dot x_1=x_2\\
     &\dot x_2=u\\
     &y=\frac{1}{2}\|x_1\|^2
 \end{align}
 where $x_1\in\mathbb{R}^n$ represents the position of the vehicle, $x_2\in\mathbb{R}^n$ is its linear velocity, and $u$ is the corresponding inertial acceleration. The output $y$ represents the (half squared) position range to the origin. The vehicle's dynamics are written as in \eqref{sys:LTI} with $x=(x_1^\top,x_2^\top)^{\top}\in\Rat{2n}$ and
\begin{align*}
&A=\left[\begin{matrix}
0&I_n\\
0&0\end{matrix}\right],\qquad 
B=\left[\begin{matrix}
0\\
I_n\end{matrix}\right],\qquad C=\left[\begin{matrix}
I_n&0\\
0&0\end{matrix}\right]
\end{align*}
 It is easy to verify that $A^2=0$ and that $C_{m}=0$ with $m=3$, i.e., $C_3=0$. Indeed, according to \eqref{Ci} we have
\begin{equation*}
\begin{array}{lll}
C_1=\left[\begin{matrix}
 0 & I_n  \\
 I_n & 0
\end{matrix}\right],&
C_2=\left[\begin{matrix}
 0 & 0  \\
 0 & 2I_n 
\end{matrix}\right],& C_3=0.
\end{array}
\end{equation*}

 Then, from \eqref{matrix:Au:ext}, we obtain the extended matrix
\begin{align*}
\mathcal{A}(u):=
\begin{bmatrix}
0&1&0&u^{\top}B^{\top}C_0\\
0&0&1&u^{\top}B^{\top}C_1\\
0&0&0&u^{\top}B^{\top}C_{2}\\
0&0&0&A
\end{bmatrix}.
\end{align*}
According to \eqref{fun:r:2} we can calculate $r_1(t)=0$, $r_2(t)=[u^\top(t),0]$, and $r_{3}(t)=[\dot{u}^{\top}(t),3u^{\top}(t)]$. Since, the matrix $A$ has real eigenvalues, a sufficient condition for the uniform observability of the extended system follows from Proposition \ref{PE:cond:4}:
$$\int_t^{t+\delta}\begin{bmatrix}
\dot{u}(s)\dot{u}^\top(s)&3\dot{u}(s)u^\top(s)\\
3{u}(s)\dot{u}^\top(s)&9u(s)u^\top(s)
\end{bmatrix}\ge \mu I_{2n},$$
which is a PE condition on the acceleration and the jerk of the vehicle.
 \begin{figure}[t!]
\centering
\includegraphics[scale=0.6]{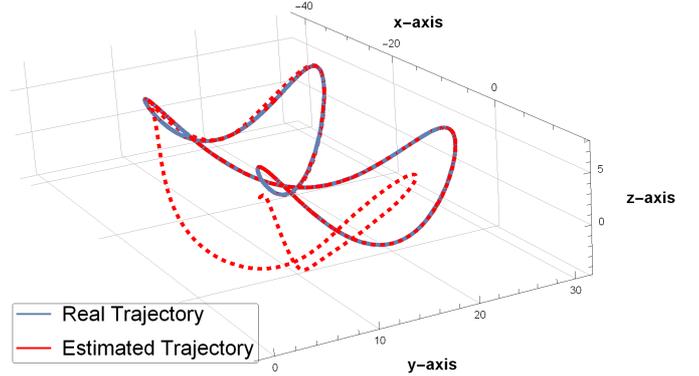}
\caption{Real and estimated trajectory of the vehicle.}
\label{fig:traj}
\end{figure}
\begin{figure}[t!]
\centering
\includegraphics[scale=0.768]{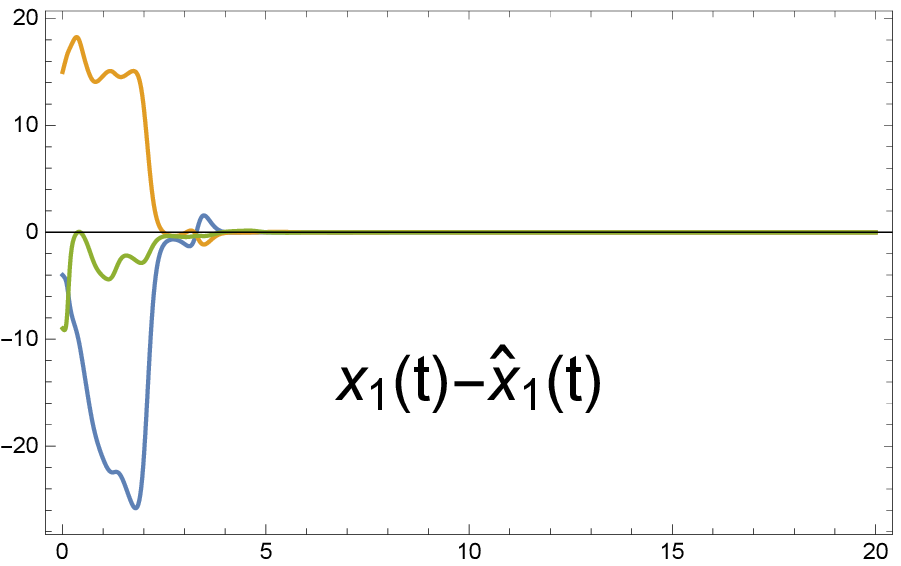}
\includegraphics[scale=0.768]{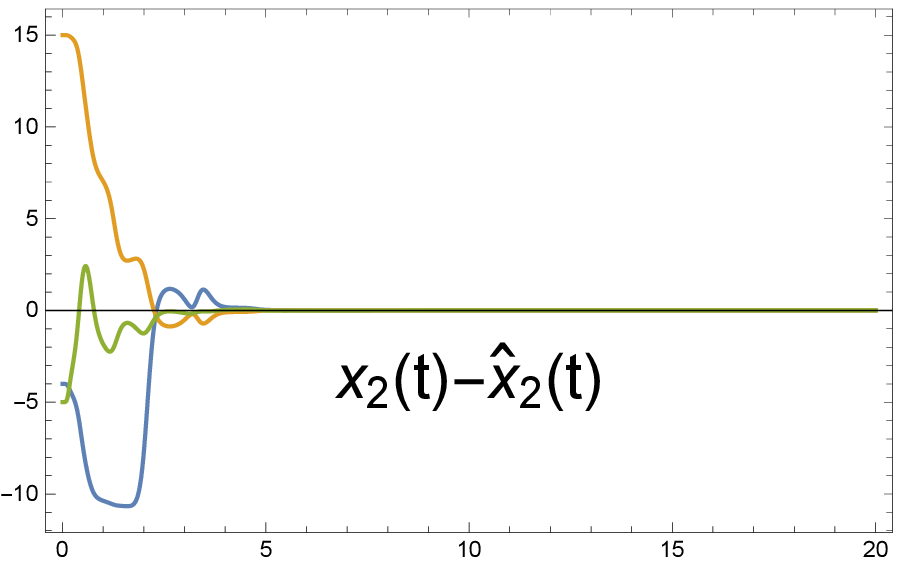}
\caption{Estimation errors for the position $x_1$ and velocity $x_2$.}
\label{fig:error}
\end{figure}
For simulation, we consider that the vehicle moves along the 3D trajectory $x_1(t)=(20 \cos (t)-20,10 \sin (2t)+20,-4 \cos (4 t))$ which is rich enough to satisfy the PE condition.  We perform the state estimation of the augmented state through a Kalman type observer given in Theorem \ref{Theorem:K:obs} with  $S_0=100I_{6}$, $V=0.0001I_6$, and $W=I_6$. In Figure \ref{fig:traj} the red trajectory generated by the observer of the system converges to the actual trajectory in blue. Figure \ref{fig:error} shows that the position and velocity estimation errors converge to zero.

 \section{Conclusion}
We proposed an immersion-type technique that transforms a class of linear systems with quadratic output to a new system with linear output  by adding a finite number of states to the original system. The class of linear systems considered is characterized by polynomial outputs under zero inputs, which encompasses for example nilpotent systems. Moreover, we derived persistence of excitation conditions for the admissible inputs that establish the uniform observability of the new system. The PE conditions are explicit and can be checked easily for a given input function. In future work we will address the problem of state estimation with multiple quadratic outputs and extend the current approach to systems for which $C_m\neq0$ for all $m\in\mathbb{N}$.

\cite{*}
\bibliography{ifacconf}
 
\appendix
\section{Proof of Binomial Expression \eqref{Ci:formula}.}

The proof of \eqref{Ci:formula} follows by induction. More specifically, for $i=1$, \eqref{Ci:formula} obviously holds. Suppose that for some $m\in\mathbb{N}$, $C_m=\sum_{r=0}^{m}\binom{m}{r}A^{\top^r}CA^{m-r}$. Then we have
	\begin{align*}
	C_{m+1}=&C_mA+A^{\top}C_m=\sum_{r=0}^{m}\binom{m}{r}A^{\top^r}CA^{m+1-r}\\&+\sum_{r=0}^{m}\binom{m}{r}A^{\top^{r+1}}CA^{m-r}\\
	=&A^{\top^{m+1}}C+\sum_{r=0}^{m-1}\binom{m}{r}A^{\top^{r+1}}CA^{m-r}\\&+\sum_{r=0}^{m}\binom{m}{r}A^{\top^r}CA^{m+1-r}\\
	=&CA^{m+1}+A^{\top^{m+1}}C+\sum_{r=1}^{m}\binom{m}{r-1}A^{\top^r}CA^{m+1-r}\\&+\sum_{r=1}^{m}\binom{m}{r}A^{\top^r}CA^{m+1 -r}\\
	=&CA^{m+1}+A^{\top^{m+1}}C+\sum_{r=1}^{m}\binom{m+1}{r}A^{\top^r}CA^{m+1-r} \\
	=&\sum_{r=0}^{m+1}\binom{m+1}{r}A^{\top^r}CA^{m+1-r}
	\end{align*}
	where in the last equality we have taken into account Pascal's identity $\binom{n}{k}+\binom{n}{k-1}=\binom{n+1}{k}$.

\end{document}